\documentstyle{article}
\begin{document}

\def\F{{\cal F}}
\def\be{\begin{equation}}
\def\ee{\end{equation}}
\def\bc{\begin{center}}
\def\ec{\end{center}}

\author{I. A. Taimanov}
\title
{Integrable geodesic flows of 
non-holonomic metrics}

\date{}
\maketitle

\bc
{\bf \S 1. Introduction}
\ec         

In the present article we show how to 
produce new examples of integrable dynamical systems
of differential geometry origin. 

This is based on a construction of a canonical 
Hamiltonian structure for the geodesic flows of 
Carnot--Carath\'eodory metrics
(\cite{G,VG}) via the Pontryagin maximum principle.
This Hamiltonian structure is 
achieved by introducing Lagrange multipliers bundles being the phase 
spaces of these Hamiltonian flows. These bundles are diffeomorphic 
to cotangent bundles but have another meaning. A transference 
to this phase space is given by a generalised Legendre transform. 

We analyse the geodesic flow of the left-invariant 
Carnot--Carath\'eodory metric on the three-dimensional Heisenberg group 
as a super-integrable Hamiltonian system (Theorem 1).
Moreover, its super-integ\-ra\-bi\-li\-ty explains the foliation of its phase
space into one- and two-dimensional invariant submanifolds
as it was pointed out in \cite{VG}.

The geodesic flows of left-invariant Carnot--Carath\'eodory metrics
on Lie groups reduce to equations on Lie algebras in the same manner
as the geodesic flows of left-invariant Riemannian metrics reduce to
the Euler equations on Lie algebras (\cite{A})
(Theorem 2). Most of these flows are integrable. 
Moreover, this reduction to equations 
on Lie algebras gives a Hamiltonian explanation for the description of such 
flows on three-dimensional Lie groups given in \cite{VG} by using
Euler--Lagrange equations (the will to explain this in terms of 
integrability was the starting point for the present work).

In comparison with the case of Riemannian 
geodesic flows there is another class
of invariant flows, corresponding to left-invariant metrics and 
right-invariant distributions. In \S 6 
we examine the simplest example of such flow on ${\cal H}^3$ 
and, in particular, 
show that this flow is integrable (Theorem 3).

In \S 7 we consider an observation related to a
Hamiltonian structure for the equations for the motion of a heavy
rigid body with a fixed point.

For completeness of explanation we discuss in \S 8
another approach to defining
``straight lines'' in non-holonomic geometry which does not arrive at
Hamiltonian systems but some ``straight line''
flows have important
mechanical meaning (for instance, the Chaplygin top (\cite{Ch})).

We also discuss some problems concerning dynamics and, in 
particular, integrability of  these systems
(see \S 9. Concluding remarks).

The present article is dedicated to D. V. Anosov on his 60th birthday.

\bc
{\bf \S 2. The geodesic flows of Carnot--Carath\'eodory metrics}
\ec

{\bf A. Carnot--Carath\'eodory metrics.}

Let $M^n$ be a smooth manifold of dimension $n$. 

A family $\F$ of $k$-dimensional subspaces of the tangent spaces to $M^n$
is called a $k$-dimensional smooth distribution if $\F_x$ is a smooth 
section of the Grassmann bundle on $M^n$. 

In the sequel we suppose that distributions are smooth.

Let $V_{\F}$ be the linear space spanned by vector fields tangent to $\F$.
Denote by $A_{\F}$ the algebra generated by fields from $V_{\F}$ via a 
commutation. A distribution is called 
{\it non-holonomic} if $A_{\F}$ does not
coincide with $V_{\F}$ as a linear space. Otherwise, a distribution is called 
holonomic and, by the Frobenius theorem, locally looks like the family of 
spaces tangent to the leaves of a foliation. It is easily seen that near a 
generic point the distribution 
corresponding to $A_{\F}$ is holonomic. The distribution 
$\F$ is called {\it completely non-holonomic} if the algebra $A_{\F}$      
coincides with the whole algebra of vector fields on $M^n$.
In \cite{St} it is said that such distributions satisfy the {\it bracket
generating hypothesis}.

In the sequel we assume that distributions are  
completely non-holonomic. This assumption is not very strong because
otherwise we may restrict geodesic flows to the leaves of the foliation and 
consider the restricted distributions as completely non-holonomic. 
Really, we use this for defining a Carnot--Carath\'eodory 
metric as a correct intrinsic metric only.

Thus we assume now that $M^n$ is endowed with a completely non-holonomic 
distribution $\F$. We also assume that $M^n$ is a complete Riemannian 
manifold with the metric ${\tilde g}^{ij}$.

A piecewise smooth  
curve in $M^n$ is called admissible if it is tangent to $\F$.
A  {\it Carnot--Carath\'eodory metric} $d_{CC}(x,y)$ is defined 
as follows.
Denote by $\Omega_{x,y}$ the set of admissible curves 
with ends at the points $x$ and $y$ in $M^n$. Then 
\be
d_{CC}(x,y) = \inf_{\gamma \in \Omega_{x,y}} 
{\mbox {length}}(\gamma)
\label{2.1}
\ee
with the lengths of curves taken with respect to the metric ${\tilde g}^{ij}$. 

By the Chow--Rashevskii theorem, any pair of points in 
a complete Riemannian manifold
endowed with a completely non-holonomic distribution is
connected by an admissible curve and, thus, (\ref{2.1}) correctly defines
an inner metric on $M^n$.

Carnot--Carath\'eodory metrics are simplest examples of non-holonomic metrics
which are defined (\ref{2.1}) for different choices
of $\Omega_{x,y}$ corresponding to non-integrable constraints. 
In the case of Carnot--Carath\'eodory metrics these constraints are linear
in velocities.

{\bf B. The geodesic flow of a Carnot--Carath\'eodory metric.}

By definition, the lengths of 
admissible curves depend only on the restrictions of the 
metric ${\tilde g}^{ij}$ onto $\F_x$. Denote these restricted forms by $Q_x$.
This family of bi-linear forms on $\F$ enables us to define 
the canonical mapping 
\be
g(x) : T^*M^n \rightarrow \F_x \subset TM^n
\label{2.2}
\ee
taking for $g(x)\xi \in \F_x$ a vector  
determined uniquely by the condition
\be
Q_x(Y, g(x)\xi) = \langle Y, \xi \rangle 
\ \ \ {\mbox{for every}}\ \  Y \in \F_x.
\label{2.3}
\ee  
The symmetric 
tensor $g^{ij}$ is called a {\it Carnot--Carath\'eodory metric tensor}.
It generalises a Riemannian metric tensor into which it degenerates when 
$\F_x = T_xM^n$.

A curve $\tilde \gamma$ in $T^*M^n$ is called a cotangent lift of a curve 
$\gamma$ in $M^n$ if 
\be
g(\gamma(t))\xi(t) = \frac{d \gamma(t)}{dt}
\label{2.4}
\ee
where ${\tilde\gamma}(t) = (\gamma(t), \xi(t)), \xi(t) \in 
T^*_{\gamma(t)}M^n$.

By (\ref{2.3}), in the terms of cotangent lifts 
the length of an 
admissible curve $\gamma(t)$ in $M^n$ is expressed by 
\be
L(\gamma) = \int \sqrt
{\langle g(\gamma(t))\xi(t), \xi(t) \rangle} dt
\label{2.5}
\ee
and the energy of $\gamma$ equals
\be
E(\gamma) = \frac{1}{2} \int \langle g(\gamma(t))\xi(t), \xi(t) \rangle dt.
\label{2.6}
\ee 

An admissible curve is called a 
{\it geodesic of the Carnot--Carath\'eodory metric $g^{ij}$} if locally it is 
an energy-minimising curve. 

The geodesics of the Carnot--Carath\'eodory metric $g^{ij}$ are described by 
the Euler--Lagrange equations for a Lagrange function
\be
L(x, \dot x) = \frac{1}{2}{\tilde g}_{pq}{\dot x}^p{\dot x}^q + 
\sum_{\alpha=1}^{n-k} \mu_{\alpha} \langle \dot x, \omega^{(\alpha)} \rangle
\label{2.7}
\ee
where $\omega^{(1)}, \dots, \omega^{(n-k)}$ is a basis for $\F_x^{\perp}$.
Here $\F_x^{\perp}$ is the annihilator of $\F_x$, i.e., the subset of 
$T^*_xM^n$ formed by covectors $\xi$ such that $\langle \xi, v \rangle = 0$
for every $v \in \F_x$.

These equations are written as 
$$
\frac{d}{dt} \frac{\partial L}{\partial {\dot x}^i} - 
\frac{\partial L}{\partial x^i} =
$$
\be
\{\frac{\partial {\tilde g}_{iq}}{\partial x^p}{\dot x}^p{\dot x}^q + 
{\tilde g}_{iq}{\ddot x}^q
+ \sum_{\alpha} ({\dot \mu}_{\alpha} \omega^{(\alpha)}_i + \mu_{\alpha}
\frac{\partial \omega^{(\alpha)}}{\partial x^p}{\dot x}^p)\} -
\label{2.8}
\ee
$$
\{\frac{1}{2}\frac{\partial {\tilde g}_{pq}}{\partial x^i}
{\dot x}^p{\dot x}^q + \sum_{\alpha} 
\mu_{\alpha}{\dot x}^p \frac{\partial \omega^{(\alpha)}_p}{\partial x^i}\} 
= 0,
$$
\be
\frac{\partial L}{\partial \mu_{\alpha}} = 
\langle {\dot x}, \omega^{(\alpha)} \rangle = 0.
\label{2.9}
\ee 

In fact, although the Riemann metric tensor ${\tilde g}^{ij}$ enters 
these equations, the geodesic flow is determined by the restriction of 
this tensor onto the distribution, the Carnot--Carath\'eodory metric tensor
$g^{ij}$, only.

\bc
{\bf \S 3. The Pontryagin maximum principle and a Hamiltonian structure for
the geodesic flows of a Carnot--Carath\'eodory metric. A generalised
Legendre transform}
\ec

{\bf A. The Pontryagin maximum principle and the geodesics of
Carnot--Carath\'eodory metrics.}

In \cite{St} Strichartz had shown that the geodesic flow of a 
Carnot-Carath\'eodory metric is described by equations derived from the
Pontryagin maximum principle.

Explain his idea in brief. First, roughly quote the
Pontryagin maximum principle in a weak form sufficient for our study 
referring to \cite{Ces} (Theorem 5.1) 
for the absolutely rigorous statement.

{\bf Pontryagin Theorem.}
{\sl Consider the minimum problem for a functional
\be
I[x(t),u(t)] = \int_{t_1}^{t_2} f^0(x,u) dt
\label{3.1}
\ee
in the class of admissible functions $(x(t),u(t))$ such that
\be
{\dot x}^k = f^k(x,u)
\label{3.2}
\ee
and some constraints $x \in A, u \in U$ hold.

Introduce the functions 
\be
{\tilde H}(x,u,{\tilde \lambda}) = \lambda_0 f^0(x,u) + 
\lambda_1 f^1(x,u) + \dots + \lambda_n f^n(x,u)
\label{3.3}
\ee
and
\be 
M(x,{\tilde \lambda}) = \inf_{u \in U} {\tilde H}(x,u,{\tilde \lambda})
\label{3.4}
\ee
where ${\tilde \lambda} = (\lambda_0, \lambda_1, \dots, \lambda_n)$.

Let $(x(t), u(t))$ be a solution to this minimum problem.
Then there exists an absolutely continuos vector function 
${\tilde \lambda}(t)$ such that

(i) $\lambda_0 = \mbox{const}, \lambda_0 \geq 0$, and
\be
\frac{d \lambda_i}{dt} = 
- \frac{\partial {\tilde H}(x(t), u(t), {\tilde\lambda}(t))}{\partial x^i} ;
\label{3.5}
\ee

(ii) ${\tilde H}(x(t), u(t), {\tilde\lambda}(t)) = 
M(x(t), {\tilde\lambda})(t)$
for every $t \in [t_1,t_2]$ ;

(iii) $M(x(t), {\tilde \lambda}(t)) = \mbox{const}$.}   

Now it suffices only to notice that
in the case of the geodesic flows of Carnot--Carath\'eodory metrics we have
the minimum problem on the set of the cotangent lifts of admissible
curves on $M^n$ where we consider the variables $\xi_i$ as the control 
functions $u_i$. In this event 
\be
f^0(x,u) = \frac{1}{2} g^{ij}u_iu_j,
\label{3.6}
\ee
\be
f^i(x,u) = g^{ij}u_j,
\label{3.7}
\ee
and
\be
{\tilde H}(x,u,{\tilde \lambda}) = \frac{\lambda_0}{2}g^{ij}u_iu_j + 
g^{ij}\lambda_iu_j.
\label{3.8}
\ee

Thus, we have
\be
M(x,{\tilde \lambda}) = 
\cases{
-\frac{1}{2\lambda_0}g^{ij}\lambda_i\lambda_j, & for $\lambda_0 \neq 0$ \cr
0, & for $\lambda_0 =0$ and $g^{ij}\lambda_j \equiv 0$ \cr
-\infty, & otherwise 
}.
\label{3.9}
\ee

{\bf B. A Hamiltonian structure for the geodesic flows of 
Carnot--Cara\-th\'eo\-dory metrics.}

Let consider the bundle $\Lambda M^n \rightarrow M^n$ 
diffeomorphic to the cotangent 
bundle $T^*M^n$ via the diffeomorphism $(x,p) \leftrightarrow (x,\lambda)$
but having another sense. We call it the {\it Lagrange multipliers bundle}
on $M^n$. As well as the cotangent bundle this bundle is endowed with 
the natural symplectic structure 
generated by the form 
\be
\Omega = \sum_{i=1}^n d\lambda_i \wedge dx^i.
\label{3.9a}
\ee

Consider now on this symplectic manifold a Hamiltonian flow 
with a Hamiltonian function 
\be
H(x,\lambda) = - \frac{1}{2}g^{ij}\lambda_i\lambda_j.
\label{3.10}
\ee

{\bf Definition.} A geodesic of a Carnot--Carath\'eodory metric is called
{\it normal} if $\lambda_0 \neq 0$. 

Otherwise, if $g^{ij}\lambda_j = 0$ then $M(x,\lambda) = 0$ 
and the Pontryagin theorem gives nothing.

{\bf Theorem HS 
(on a Hamiltonian structure for normal geodesic flow)}.
{\sl The projections of trajectories of the Hamiltonian flow on $\Lambda M^n$
with the Hamiltonian function (\ref{3.10}) are exactly the
naturally-parametrised normal geodesics of the Carnot--Carath\'eodory metric
$g^{ij}$. Moreover, $|{\dot x}|^2 = -2 H(x,\lambda)$.}

Proof of Theorem HS.

From 
(\ref{3.9}), and the homogeneity of ${\tilde H}(x,u,{\tilde \lambda})$ 
we infer that after changing of the parameter on an extremal to a multiple 
one, if that needs, we obtain an extremal ${\tilde \gamma}$ 
with $\lambda_0 = 1$.

By (\ref{3.4}), we derive
$$
\frac{\partial {\tilde H}(x,u,{\tilde \lambda})}{\partial u_i} = 0
$$
which is equivalent to
\be
g^{ij}u_j = - g^{ij}\lambda_j.
\label{3.11}
\ee
Now it follows from (\ref{3.11}) that
\be
g^{ij}u_iu_j = g^{ij}\lambda_i\lambda_j.
\label{3.12}
\ee
 
Remind that, by (\ref{2.4}), we have
\be
g^{ij}u_j = {\dot x}^i
\label{3.13}
\ee
and together with (\ref{2.3}) this implies that
\be
|{\dot x}|^2 = {\tilde g}_{ij}{\dot x}^i{\dot x}^j = g^{ij}u_iu_j = |u|^2.
\label{3.14}
\ee

Now it follows from (\ref{3.12}) and   
the statement (iii) of the Pontryagin theorem that
the extremal ${\tilde \gamma}$ is naturally-parametrised.

Since (\ref{2.4}) and (\ref{3.11}), we have
that
\be 
{\dot x}^i = \frac{\partial H(x,\lambda)}{\partial \lambda_i}
\label{3.15}
\ee
and we may regard (\ref{3.5}) as
\be
{\dot \lambda_i} = - \frac{\partial H(x,\lambda)}{\partial x^i}.
\label{3.16}
\ee   

It is easily seen that the equations (\ref{3.15}--\ref{3.16}) form a system of 
Hamilton equations on $\Lambda M^n$ for the Hamiltonian function 
(\ref{3.10}).

Thus we prove that the normal geodesics are the projections of 
trajectories of this Hamiltonian flow.

For the converse we refer to \cite{H}, where a comprehensive
examination of analytical properties of the energy functional for 
Carnot--Carath\'eodory metrics is given.

This completes the proof of the theorem.

The variables $\lambda_1, \dots, \lambda_n$ have no physical meanings in
comparison with the momenta $u_1, \dots, u_n$. 
The correspondence between them and 
velocities is given by (\ref{3.15}) 
and is one-to-one in the case when the 
form $g^{ij}$ is non-degenerate, i.e., for Riemannian metrics only. 
In this case the existence of a Hamiltonian formalism for the geodesic flows
of Riemannian metrics follows from both the Legendre transform and the 
Pontryagin theorem. Thus the  transference to the new variables 
$(x, \lambda)$ has to be regarded as a {\it generalised Legendre transform}.  

Theorem HS is contained implicitly in \cite{St}. 
However, this fact did not attract attention of specialists
on integrable systems because mostly 
the problem of regularity of geodesics and that of exponential maps were 
in study (\cite{H,Tay}).  
From the analytical point of view this also coincides with the
introducing a Hamiltonian in the ``vakonomic mechanics'' (\cite{AKN})
where the momenta are introduced by the implicit theorem procedure
and the Hamiltonian system is also regarded as a system on a cotangent 
bundle. However, the difference between $T^*M^n$ and
$\Lambda M^n$ is essential for applications 
to mechanics and physics because, at least,  
the variables $\lambda_1, \dots, \lambda_n$ are not observable.

Although for a long time it had being assumed that all geodesics are normal, 
recently Montgomery had showed that abnormal geodesics exist (\cite{Mon}).
Geodesics found by him do not admit end-point $\F$-tangent 
perturbations and,
thus, are solutions to any variation problem on the space of 
admissible curves.

But if the space ${\cal V}_{\F} + [{\cal V}_{\F}, {\cal V}_{\F}]$ 
coincides with the whole algebra of vector fields on $M^n$ then every
geodesic is normal (\cite{St}).

In the sequel speaking about geodesic flows we shall mean by them
normal geodesic flows.

\bc
{\bf \S 4. Integrability of the geodesic flow of the 
left-invariant Carnot--Carath\'eodory metric on the three-dimensional
Heisenberg group}
\ec 

We mean by a left-invariant Carnot--Carath\'eodory metric 
a left-invariant metric restricted onto a left-invariant distribution.
It is known that such metric on the three-dimensional Heisenberg group is 
unique up to isomorphism (\cite{VG}). 

This flow is described in \cite{B2,H,VG,VoG}, somewhere  
with its generalizations.
However, it is nowhere regarded as a completely integrable Hamiltonian system.

The three-dimensional Heisenberg group ${\cal H}^3$ 
is the group of matrices
\be
\left(
\begin{array}{ccc}
1 & x & z \\
0 & 1 & y \\
0 & 0 & 1
\end{array}
\right)
\label{4.1}
\ee
with respect to a multiplication, where $x, y, z \in {\bf R}$.
Its Lie algebra ${\cal L}$ is spanned by the following elements
\be
e_1 = 
\left(
\begin{array}{ccc}
0 & 1 & 0 \\
0 & 0 & 0 \\
0 & 0 & 0
\end{array}
\right), \
e_2 =
\left(
\begin{array}{ccc}
0 & 0 & 0 \\
0 & 0 & 1 \\
0 & 0 & 0 
\end{array}
\right), \
e_3 =
\left(
\begin{array}{ccc}
0 & 0 & 1 \\
0 & 0 & 0 \\
0 & 0 & 0
\end{array}
\right).
\label{4.2}
\ee
Denote by ${\cal L}_0$ the linear subspace spanned by $e_1$ and $e_2$.

The group ${\cal H}^3$ acts on itself by 
the left and right translations :
$$
L_g : {\cal H}^3 \rightarrow {\cal H}^3 : \ \ L_g(h) = gh,
$$
$$
R_g : {\cal H}^3 \rightarrow {\cal H}^3 : \ \ R_g(h) = hg.
$$
The left-invariant distribution generated by ${\cal L}_0$ consists of
the 2-planes ${\cal F}_x = L_{g*} {\cal L}_0$. 

Since the following commutation relations
$$
[e_1,e_2] = e_3, \ \ 
[e_1,e_3] = [e_2,e_3] = 0
$$ 
hold, this distribution is completely non-holonomic.

Consider the left-invariant metric on ${\cal H}^3$ which at the unit of 
the group takes the form 
\be
(e_i, e_j) = \delta_{ij}.
\label{4.3}
\ee
Identify ${\cal H}^3$ with ${\bf R}^3$
by the diffeomorphism which assigns to the matrix (\ref{4.1})
the point in ${\bf R}^3$ with the coordinates $(x,y,z)$. 
Thus we identify the tangent space at every point of ${\cal H}^3$ with
the vector space generated by the matrices (\ref{4.2}).
In this event the left translations act on $T{\cal H}^3$ as follows :
\be 
L_{g*}(e_1) = e_1, \ \ 
L_{g*}(e_2) = e_2 + xe_3, \
L_{g*}(e_3) = e_3
\label{4.4}
\ee
where $g$ is the element of ${\cal H}^3$ given by the matrix (\ref{4.1}).
It follows from (\ref{4.4}) that in these coordinates the mapping
$L_{g*} : T_e{\cal H}^3 \rightarrow T_g{\cal H}^3$ is written as
\be
L_{g*} =
\left(
\begin{array}{ccc}
1 & 0 & 0 \\
0 & 1 & 0 \\
0 & x & 1
\end{array}
\right)
\label{4.4a}
\ee
where $e$ is the unit of ${\cal H}^3$.
Thus we have 
\be
[g^{ij}(x,y,z)] = (L_{g*})^{**} \cdot [g^{ij}(0,0,0)] \cdot
(L_{g*})^*
\label{4.4b}
\ee
and, since we put
\be 
g^{ij}(0,0,0) =
\left(
\begin{array}{ccc}
1 & 0 & 0 \\
0 & 1 & 0 \\
0 & 0 & 0
\end{array}
\right),
\label{4.4c}
\ee
we derive
\be
[g^{ij}(x,y,z)] = 
\left(
\begin{array}{ccc}
1 & 0 & 0 \\
0 & 1 & x \\
0 & x & x^2
\end{array}
\right).
\label{4.5}
\ee

The left-invariant Riemannian metric on ${\cal H}^3$ in these coordinates
takes the form
\be
[{\tilde g}_{ij}(x,y,z)] =
\left(
\begin{array}{ccc}
1 & 0 & 0 \\
0 & (1+x^2) & -x \\
0 & -x & 1
\end{array}
\right).
\label{4.5a}
\ee

Now, it follows from Theorem HS that 
the geodesic flow of the left-invariant Carnot--Carath\'eodory metric 
corresponding to the Riemannian metric (\ref{4.5a}) and the distribution
$L_{g*}{\cal L}_0$ is Hamiltonian on $\Lambda {\cal H}^3$ 
with the following Hamiltonian function
\be
H(q,\lambda) = \frac{1}{2}(\lambda_1^2 + \lambda_2^2 + x^2\lambda_3^2 +
2x\lambda_2\lambda_3)
\label{4.6}
\ee
where
$q = (x,y,z)$.

The Hamilton  equations for ({\ref{4.6}) look simply :
$$
{\dot x} = \frac{\partial H}{\partial \lambda_1} = \lambda_1, \ \ 
{\dot y} = \frac{\partial H}{\partial \lambda_2} = 
\lambda_2 + x\lambda_3,
$$
\be 
{\dot z} = \frac{\partial H}{\partial \lambda_3} =
x\lambda_2 + x^2\lambda_3,\ \  
{\dot \lambda}_1 = - \frac{\partial H}{\partial x} = 
- x\lambda_3^2 - \lambda_2\lambda_3, 
\label{4.7}
\ee
$$
{\dot \lambda}_2 =  - \frac{\partial H}{\partial y} = 0,\ \
{\dot \lambda}_3 = -\frac{\partial H}{\partial z} = 0,
$$
and we immediately infer from (\ref{4.7}) that this Hamiltonian system is 
completely integrable because it has three first integrals 
\be
I_1=H,  \ \ I_2 = \lambda_2, \ \ I_3=\lambda_3
\label{4.7a}
\ee
which are in involution and
functionally independent almost everywhere. In particular, these integrals
are functionally independent in the domain $H \neq 0$.

Moreover, 
since, by (\ref{4.7}), 
\be
\lambda_3 = const, \ \ 
{\dot z} =  x{\dot y}
\label{4.8}
\ee
along trajectories of the flow, we may restrict this flow onto the level set 
$\{\lambda_3 = C = \mbox{const}\}$ and project this restriction of the flow
onto the plane $(x,y)$. Denote this system by ${\cal P}_C$ and notice that
it is  defined on the $4$-dimensional
symplectic manifold ${\cal M}_C$ 
diffeomorphic to the cotangent bundle to the $2$-plane
with the coordinates $(x,y)$ but with another Poisson structure. 

Introduce the new variables on ${\cal M}_C$
\be
u = \lambda_1, \ \ v = \lambda_2 + x \lambda_3.
\label{4.9}
\ee
Then, by (\ref{3.9a}) and (\ref{4.9}), the Poisson structure 
on ${\cal M}_C$ induced from $\Lambda {\cal H}^3$ is written as
\be 
\{x,u\} = \{y,v\} = 1, \{u,v\} = - C \ (= - \lambda_3),
\label{4.10}
\ee 
$$
\{x,v\} = \{y,u\} = \{x,y\} = 0
$$
in the coordinates $(x,y,u,v)$.
The flow ${\cal P}_C$ is also a Hamiltonian 
system with the following Hamiltonian function
\be
H(x,y,u,v) = \frac{u^2+v^2}{2}.
\label{4.11}
\ee

It is easily seen that the flow (\ref{4.10}--\ref{4.11}) describes 
nothing else but the motion of a charged particle on the Euclidean 
plane $(x,y)$ in the constant magnetic field $F = -\lambda_3 dx \wedge dy$
(\cite{NS}). This system has three first integrals which are functionally 
independent almost everywhere 
\be
{\tilde I}_1 = H, \ \ {\tilde I}_2 = \lambda_3 x - v, \ \ 
{\tilde I}_3 = \lambda_3 y + u,
\label{4.12}
\ee
and, thus, we conclude that these systems are 
super-integrable, i.e. have more
first integrals than $\dim {\cal M}_C/2$. This is also true for the main flow.

Hence, we conclude

{\bf Theorem 1.}
{\sl 1) The geodesic flow of the left-invariant Carnot--Carath\'eodory
metric, on ${\cal H}^3$, corresponding to the left-invariant Riemannian 
metric (\ref{4.5a}) and the left-invariant distribution $L_{g*}{\cal L}_0$
is a Hamiltonian system on $\Lambda {\cal H}^3$ which
is integrable in the Liouville sense via the 
first integrals (\ref{4.7a}) which are in involution and 
functionally independent almost everywhere.

Moreover, since this flow possesses the fourth first integral 
$I_4 = \lambda_3y + \lambda_1$ functionally independent on 
(\ref{4.7a}), it is super-integrable and the subset $\{ \lambda_3 
\neq 0 \}$ of its phase space is foliated into 
$2$-dimensional invariant Liouville tori
$S^1 \times {\bf R}$.

2) Let restrict this geodesic flow
onto the level set 
$\{\lambda_3 = C\}$ and project this restriction of the flow onto
the plane $(x,y)$. The flow ${\cal P}_C$ constructed by this procedure
is equivalent to the Hamiltonian system 
describing the motion of a charged particle  on the Euclidean $2$-plane 
$(x,y)$ in the constant magnetic field $F = -\lambda_3 dx \wedge dy$. 
This flow is super-integrable and for $\lambda_3 \neq 0$
its phase space ${\cal M}_C$ is foliated into closed trajectories.}

\bc
{\bf \S 5. Left-invariant Carnot--Carath\'eodory metrics on Lie groups
and their geodesic flows} 
\ec

Let ${\cal G}$ be a Lie group, let $G$ be its Lie algebra, and let 
$G_0$ be a subspace, of $G$, generating $G$. Take a scalar product
${\cal J}$ in $G$ and decompose $G$ into a direct sum of $G_0$ and 
its orthogonal complement :
$$
G = G_0 \oplus G_0^{\perp}.
$$  
Take  another bi-linear from ${\cal J}_0$ on $G$ uniquely defined by 
the following conditions 
$$
{\cal J}_0(x,y) = 0\ \  \mbox{for every} \ \ x \in G_0^{\perp}, y \in G,
$$
$$
{\cal J}_0(x,y) = {\cal J}(x,y) \ \ \mbox{for every} \ \ x, y \in G_0.
$$

Now, take the left-invariant distribution $L_{g*}G_0$ on ${\cal G}$ and
the left-invariant Riemannian metric generated by ${\cal J}$. To this pair
there is uniquely assigned a left-invariant Carnot--Carath\'eodory metric.

We will not plunge into the details but only mention that for they 
same reasoning as the geodesic flow of a left-invariant Riemannian metric
on a Lie group reduces to equations on its Lie co-algebra, 
Euler equations (\cite{A}, 
the following theorem holds.

{\bf Theorem 2.}
{\sl The geodesic flow of a left-invariant 
Carnot--Carath\'eodory metric reduces to 
the following equations on $G^*$
\be
{\dot M} = ad^*_{\omega} M
\label{11.1}
\ee
where $\omega = L_{g^{-1}*} {\dot g} \in G$ and $M = {\cal J}_0 \omega$.
Here ${\cal J}_0$ is regarded as an operator ${\cal J}_0 : G \rightarrow G^*$
acting as ${\cal J}_0(x,y) = \langle {\cal J}_0x, y \rangle $.

If there exists an invariant non-degenerate bi-linear form on $G$ 
we may identify  $G$ and $G^*$ and the equations (\ref{11.1}) take the form}
\be
{\dot M} = [\omega, M], \ \ M={\cal J}_0 \omega.
\label{11.2}
\ee

Since these flows possess commutation representations (\ref{11.2}), they
give a lot of new examples of integrable Hamiltonian systems and 
the well-known methods of
integrating Euler equations on Lie algebras are
immediately generalised for them. 

Consider the simplest example. Let ${\cal G} = SO(3)$, 
let $e_1, e_2, e_3$ be its generators satisfying the commutation relations
$$
[e_i, e_j] = \varepsilon_{ijk} e_k,
$$
let ${\cal J} = \mbox{diag}(1,1,1)$, and let $G_0$ be spanned by 
$e_1$ and $e_2$.
In this case, the equations (\ref{11.2}) 
has two first integrals $\langle {\cal J} x,x \rangle$ and 
$\langle {\cal J}_0 x, x \rangle$ and, thus, are completely integrable.

We would like to notice that 
passing from the Lagrange equations (\ref{2.8}--\ref{2.9}) to the Hamiltonian 
equations (\ref{11.1}) simplifies the research of left-invariant
Carnot--Carath\'eo\-do\-ry geodesic flows and explains 
the integrable behaviour
of such flows on three-dimen\-sio\-nal Lie algebras given in \cite{VG}.

\bc
{\bf \S 6. The geodesic flow of the Carnot-Carath\'eodory metric on
${\cal H}^3$ corresponding to a left-invariant metric and 
a right-invariant distribution}
\ec

Since a Carnot--Carath\'edory metric is defined by a pair of 
objects (a metric and a distribution, there are other classes of invariant 
Carnot--Carath\'edory metrics corresponding to metrics and distributions 
invariant with respect to different actions of a Lie group ${\cal G}$. 

Examine, for instance, the 
geodesic flow of the
Carnot--Carath\'eodory metric on ${\cal H}^3$ which corresponds to 
the left-invariant metric (\ref{4.5a}) and the right-invariant distribution 
$R_{g*}{\cal L}_0$ and show that it is integrable in the Liouville sense.

By (\ref{2.3}) and (\ref{4.5a}), we have
$$
g^{11} = (1+x^2) V, \ \ 
g^{12} = g^{21} = xyV,
$$
\be
g^{22} = (1+y^2)V, \ \
g^{13} = g^{31} = y(1+x^2)V, \ \
\label{5.2}
\ee
$$
g^{23} = g^{32} = xy^2V,\ \
g^{33} = y^2(1+x^2)V, \ \
V = \frac{1}{1+x^2+y^2}.
$$

This system has two obvious first integrals : $I_1=H$ and $I_2=\lambda_3$.
As in \S 4, restrict this flow onto the level set 
$\{\lambda_3 = C = \mbox{const}\}$ and successively project this restriction 
of the flow onto the plane $(x,y)$. Thus we obtain a Hamiltonian system 
${\cal R}_C$ on the 
$4$-dimensional symplectic manifold ${\cal M}_C$.
Introduce the new variables 
$$
u = \lambda_1 + y \lambda_3,\ \  v = \lambda_2.
$$
Then the Poisson structure on ${\cal M}_C$ is given by
\be
\{x,u\} = \{y,v\} = 1, \ \ \{u,v\} = C \ (= \lambda_3),
\label{5.3}
\ee
$$
\{x,v\} = \{y,u\} = \{x,y\} = 0.
$$
The Hamiltonian functions are written as
\be
H(x,y,u,v) = \frac{1}{2(1+x^2+y^2)}((1+x^2)u^2 + 2xyuv +(1+y^2)v^2).
\label{5.4}
\ee

By the same reasoning as in the proof of Theorem 2 we conclude that
this flow is equivalent to a Hamiltonian system describing the motion
of a charged particle on the $2$-plane with the Riemannian metric
\be
(1+y^2)dx^2 - 2 xydxdy + (1+x^2)dy^2
\label{5.5}
\ee
in the constant magnetic field $F = \lambda_3 dx \wedge dy$.

In the polar coordinates $(r,\varphi)$, where $x = r\cos\varphi$ and
$y = r\sin\varphi$, the metric (\ref{5.5}) is written as
\be
dr^2 + (r^2 + r^4) d\varphi^2
\label{5.6}
\ee
and we infer that the flow ${\cal R}_C$ is Hamiltonian,
$$
\frac{df}{dt} = \{f,H\},
$$
with
the  Hamiltonian function
\be
H(r,\varphi,p_r,p_{\varphi}) =
\frac{1}{2}(p_r^2 + \frac{p_{\varphi}^2}{r^2 + r^4}) 
\label{5.7}
\ee
and the following Poisson structure on ${\cal M}_C$
\be
\{r,p_r\} = \{\varphi,p_{\varphi}\} = 1, \ \ \{p_r,p_{\varphi}\} = C,
\label{5.8}
\ee
$$
\{r,p_{\varphi}\} = \{\varphi,p_r\} = \{r,\varphi\} = 0.
$$
This flow is defined on the four-dimensional symplectic manifold ${\cal M}_C$
and has two functionally independent first integrals 
${\tilde I}_1 = H$ and ${\tilde I}_2 = p_{\varphi} + Cr^2/2$.
It is clear that these functions are also first integrals of the main
geodesic flow. In the initial coordinates the integral 
${\tilde I}_2$ takes the form
$$
{\tilde I}_2  = C\frac{x^2+y^2}{2} + xv - yu.
$$

We conclude

{\bf Theorem 3.}
{\sl 1) The geodesic flow, the Carnot--Carath\'eodory metric corresponding
to the left-invariant Riemannian metric (\ref{4.5a}) and the right-invariant
distribution $R_{g*}{\cal L}_0$, is a Hamiltonian system on $\Lambda 
{\cal H}^3$ with the following Hamiltonian function
\be
H(q,\lambda) = 
\frac{1}{2(1+x^2+y^2)}((1+x^2)\lambda_1^2 + 
(1+y^2)\lambda_2^2 + 
\label{5,9}
\ee
$$
y^2(1+x^2)\lambda_3^2 + 2xy\lambda_1\lambda_2 +
2y(1+x^2)\lambda_1\lambda_3 + 2xy^2\lambda_2\lambda_3)
$$
where $q = (x,y,z)$ ;

2) This flow possesses the first integrals 
$$
I_1 = H, \ \ I_2 = \lambda_3, \ \ I_3 = \lambda_3\frac{x^2-y^2}{2} + 
x\lambda_2 - y\lambda_1
$$ 
which are involutive and functionally independent almost everywhere.
Hence, this flow is integrable in the Liouville sense ;

2) Let restrict this flow onto the level set $\{ \lambda_3 = \mbox{const}\}$
and project the restriction of the flow onto the plane $(x,y)$. The flow
${\cal R}_C$ constructed by this procedure is equivalent to a Hamiltonian 
system describing the motion of a charged particle on the $2$-plane with 
the Riemannian metric (\ref{5.5}) in the constant magnetic field
$\lambda_3 dx \wedge dy$.} 

\bc
{\bf \S 7. The equations for the motion of a heavy rigid body with
a fixed point}
\ec

First, remind that the Lie algebra $e(3)$ of the group of motions of the 
three-dimensional Euclidean space, $E(3)$, is spanned by the elements
$e_1, e_2, e_3, f_1, f_2$, and $f_3$ meeting the following commutation 
relations
\be
[e_i, e_j] = \varepsilon_{ijk} e_k, \ \
[e_i, f_j] = \varepsilon_{ijk} f_k, \ \
[f_i, f_j] = 0.
\label{19.1}
\ee
Denote by $m_i$ and $\gamma_j$ the adjoint basis in the co-algebra
$e^*(3)$. The relations (\ref{19.1}) determine the Lie--Poisson structure
on the space of functions on $e(3)$ as follows
\be
\{m_i, m_j\} = \varepsilon_{ijk} m_k, \ \
\{m_i, \gamma_j\} = \varepsilon_{ijk} \gamma_k, \ \
\{\gamma_i, \gamma_j\} = 0.
\label{19.2}
\ee

As in the case of the geodesic flows of left-invariant metrics on Lie groups,
any Lagrangian system corresponding to a left-invariant metric and a
``left-invariant'' potential field on the group $E(3)$ reduces to
a Hamiltonian system 
$$
\frac{df}{dt} = \{f, H\}
$$
on the algebra $e(3)$ with a Hamiltonian
\be
H(m,\gamma) = 
\frac{1}{2} a^{ij}m_im_j + \frac{1}{2} b^{ij}(m_i \gamma_j + m_j \gamma_i)
+ \frac{1}{2} c^{ij} \gamma_i \gamma_j + V(m, \gamma)   
\label{19.3}
\ee
where the matrices $a^{ij}, b^{ij}$, and $c^{ij}$ are symmetric and
$V$ is a linear function in $m$ and $\gamma$. Here we call a potential 
field $V(q)$ ``left-invariant'' if its gradient is left-invariant.
Thus, denoting by $x^i$ and $y^j$ 
the local coordinates corresponding to $e_i$ and $f_j$, we conclude that 
$$
V(m, \gamma) = 
\sum_{i=1}^3 
\left( \frac{\partial V(0)}{\partial x^i} m_i
+ \frac{\partial V(0)}{\partial y^j} \gamma_j
\right)
$$ 
and a ``left-invariant'' potential field is determined uniquely 
by its gradient at the unit of the group.

The Kirhgoff equations for the free motion of a rigid body in a liquid
correspond to $V \equiv 0$ (\cite{NS}). In this case the configuration 
space is the whole group $E(3)$.  

By Theorem HS, in the Hamiltonian formalism constraints contribute to 
Hamiltonian equations via a Hamiltonian function. Moreover, Theorem HS
also holds for holonomic constraints.

Hence, starting from the problem of the free motion of 
a heavy body in a potential field with the configuration space $E(3)$ we
pose holonomic constraints by fixing a point of the body. In this case the 
configuration space is homeomorphic to $SO(3)$ but 
the Euler equations (see \S 5) are still written for the algebra $e(3)$
and correspond to the Hamiltonian function
\be
H(m, \gamma) = \frac{(Im,m)}{2} + (r,\gamma).
\label{19.4}
\ee

Thus, we obtain the well-known Hamiltonian formalism for 
this problem. 

First, it was derived in physical terms. 
Here $I$ is the inertia tensor, 
$m$ is the angular momentum, 
$\gamma$ is the vector in the direction of
gravity, and $r$ is the centre of mass. 
All coordinates are taken with respect to the orthogonal 
frame attached to the body with the fixed point the centre of coordinates.

\bc
{\bf \S 8. Another definition of ``straight lines'' and problems of mechanics}
\ec

There is another way to define  
``straight lines'' in non-holonomic geometry. We discuss only the case of
constraints linear in velocities.

Roughly speaking, geodesics of Carnot--Carath\'eodory metric are
solutions of the following variation problem. Let $L(x,{\dot x}) = 
{\tilde g}_{ij}{\dot x}^i{\dot x}^j dt$ be an energy functional on a 
suitable space of curves (periodic or with fixed end-points) 
in manifold $M^n$ and ${\tilde g}_{ij}$ be a 
Riemannian metric tensor. A geodesic $\gamma(t)$ of 
the Carnot--Carath\'eodory metric
corresponding to ${\tilde g}_{ij}$ and a distribution $\F$ are 
local extremals of this functional with respect to 
the set of variations of the form $\gamma(t,u)$
where $\gamma(t,o) = \gamma(t)$ and 
\be
\frac{\partial \gamma}{\partial t}(t,u) \in \F_{\gamma(t,u)}.
\label{12.1}
\ee
Of course, all variations belong to the space of curves in study.

But the condition (\ref{12.1}) can be changed to another one :
\be
\frac{\partial \gamma}{\partial u}(t,0) \in \F_{\gamma(t,o)}.
\label{12.2}
\ee 

In this case the equations for ``straight lines'' are written as
\be
\frac{d}{dt} \frac{\partial L}{\partial {\dot x}^i} - 
\frac{\partial L}{\partial x^i} = \sum_{\alpha} \mu_{\alpha} 
\omega^{(\alpha)}_i
\label{12.3}
\ee
where
\be
\langle \omega^{(\alpha)}, {\dot x} \rangle = 
\omega^{(\alpha)}_i {\dot x}^i = 0
\label{12.4}
\ee
are the constraints and the coefficients $\mu_{\alpha}$ are derived from
the condition that the relations (\ref{12.4}) hold.
One can easily see the difference of these equations from the equations 
(\ref{2.8}--\ref{2.9}).

Notice also that these ``straight lines'' are determined by the Lagrange 
function defined on the whole tangent space $TM^n$ but not only on
$\F$.

The most famous example of such system is the Chaplygin top, the 
dynamically-asymmetric ball rolling on the horizontal plane and with the
center of mass coinciding with its geometric centre (\cite{Ch}).

This system was integrated by Chaplygin by no means of the Liouville
integrability theorem. Its algebraic origin was clarified by Veselov and
Veselova who regarded it such a flow of ``straight lines'' on 
$E(3)$, the Lie group of the motions of the three-dimensional Euclidean space,
endowed with a left-invariant metric and a right-invariant non-holonomic
distribution (\cite{VV}). They also generalise this system on arbitrary Lie
groups and succeed in generalising the Chaplygin integration method 
for three-dimensional groups.

\bc
{\bf \$ 9. Concluding remarks}
\ec

1) The procedure of constructing a Hamiltonian structure is generalised
for non-holonomic systems with Lagrange functions of the form
$$
L(x,{\dot x}) = \frac{|{\dot x}|^2}{2} + U(x)
$$
in the usual manner. We already mentioned this in \S 7. 

2) We did not find yet considerably new examples of manifolds which
admit the integrable geodesic flow of a Carnot--Carath\'eodory metric 
and have no Riemannian metrics with integrable geodesic flows. Nevertheless,
we would like to mention that all methods of finding topological 
obstructions to integrability of the geodesic flows of Riemannian metrics
(\cite{Pat,T1}, see also \cite{T2})
fail in the case of Carnot--Carath\'eodory metrics. 
The reason for this is clear. These methods use a compactness of the level
set of a Hamiltonian function which can not be compact for 
Carnot--Carath\'eodory metrics. Thus one can expect that the class of 
manifolds admitting integrable Carnot--Carath\'eodory geodesic flows
is wider than the class of manifolds admitting 
integrable Riemannian geodesic flows.

The lack of compactness of the level sets of a Hamiltonian function also 
obstructs us to define entropy characteristics of such flows in the usual 
manner. 

However, in \S 5 we give an example of an integrable flow on $SO(3)$
with compact level sets of the first integral
$({\cal J} x,x)$ which is not a Hamiltonian function. 
This situation is typical for the geodesic flows
of left-invariant metrics on compact Lie groups but it is not generic
as one can see from the example given in \S 6. 

3) Consider the simplest example of 
degeneration of integrable Riemannian geodesic flows 
into an integrable Carnot--Carath\'eodory geodesic flow.

Take the Lie group $SO(3)$ and denote by $e_1,e_2$, and 
$e_3$ the generators of its Lie algebra $so(3)$
with the following commutation relations
$$
[e_i, e_j] = \varepsilon_{ijk}e_k. 
$$
Denote by ${\tilde {\cal L}}_0$ the subspace of $so(3)$ spanned by 
$e_1$ and $e_2$.
Consider the family of left-invariant metrics generated by the metrics on
$so(3)$ of the form 
$$
G_D = 
\left(
\begin{array}{ccc}
1 & 0 & 0 \\
0 & 1 & 0 \\
0 & 0 & D
\end{array}
\right).
$$
The geodesic flows of these metrics are integrable.
Tending  $D$ to infinity, $D \rightarrow \infty$, 
we arrive at the geodesic flow of 
the Carnot--Caratheodory metric on $SO(3)$ corresponding to 
the Riemannian metric $G_1$ and the left-invariant distribution
$L_{g*}{\tilde {\cal L}}_0$ (\cite{G,M}). It is showed in \S 5 that this flow 
is integrable.

\vskip5mm

{\bf Acknowledgements.}

The author thanks I. Schmelzer 
\footnote[1]{He is also known as I. Shmel'tser 
because of the double translation from 
German to Russian and from Russian to English, see \cite{NS}.}
and A. P. Veselov for helpful conversations.  

This work was done during author's stay at the
Institute of Theoretical Physics of Freie--Universit\"at in Berlin and
was supported by the Volkswagenwerkstiftung Foundation.

\vskip5mm

Institute of Mathematics,

630090 Novosibirsk, Russia

e-mail : taimanov@math.nsc.ru


\begin{thebibliography}{999}


\bibitem{A}
V. I. Arnol'd,
Mathematical methods of classical mechanics,
Springer-Verlag : Berlin-Heidelberg-New York, 1978.

\bibitem{AKN}
V. I. Arnol'd, V. V. Kozlov, and A. I. Nejshtadt,
Mathematical aspects of classical and celestial mechanics,
in: Dynamical Systems III,
Encyclopaedia of Mathematical Sciences {\sl 3}, Springer-Verlag:
Berlin-Heidelberg-New York, 1988.

\bibitem{B1}
V. N. Berestovskii,
Homogeneous spaces with intrinsic metric,
Soviet Math. Dokl. {\bf 38} (1989), 60--63.

\bibitem{B2}
V. N. Berestovskii,
Geodesics of nonholonomic left-invariant intrinsic metrics on the Heisenberg
group and isoperimetric curves on the Minkowskii plane,
Siberian Math. Journal {\bf 35} (1994), 1--8.

\bibitem{Ces}
L. Cesari, 
Optimization theory and applications. Springer-Verlag, New York, 1983.

\bibitem{Ch}
S. A. Chaplygin,
Investigations in the dynamics of nonholonomic systems,
Gostekhizdat : Moscow-Leningrad, 1947. (Russian)

\bibitem{G}
M. Gromov,
Carnot--Carath\'eodory spaces seeing from within,
Preprint, IHES, 1994.

\bibitem{H}
U. Hamenst\"adt,
Some regularity theorems for Carnot--Carath\'eodory metrics,
J. Diff. Geometry {\bf 32} (1990), 819--850.

\bibitem{M}
J. Mitchell,
On Carnot--Carath\'eodory metrics,
J. Diff. Geometry {\bf 21} (1985), 35--45.

\bibitem{Mon}
R. Montgomery,
Abnormal minimizers,
SIAM J. Control and Optimization {\bf 32} (1994), 1605--1620.

\bibitem{NS}
S. P. Novikov and I. Shmel'tser,
Periodic solutions of the Kirhgoff equations for the free motion of
a rigid body in a liquid, and the extended 
Lyusternik--Schnirelmann--Morse (LSM) theory. I,
Functional Analysis and its Appl. {\bf 15} (1981), 197--207.

\bibitem{Pat}
G. Paternain,
On the topology of manifolds with completely integrable geodesic flows,
Ergodic Theory and Dynamical Systems {\bf 12} (1992), 109--121.

\bibitem{St}
R. S. Strichartz,
Sub-Riemannian geometry,
J. Diff. Geometry {\bf 24} (1986), 221--263 ;
J. Diff. Geometry {\bf 30} (1989), 595--596.

\bibitem{T1}
I. A. Taimanov,
Topological obstructions to integrability of geodesic flows on 
non-simply-connected manifolds,
Math. USSR Izvestiya {\bf 30} (1988), 403--409.

\bibitem{T2}
I. A. Taimanov,
The topology of Riemannian manifolds with integrable geodesic flows,
Proc. of the Steklov Institute of Mathematics {\bf 205} (1995), 139--150.

\bibitem{Tay}
T. J. S. Taylor,
Some aspects of differential geometry associated with hypoelliptic
second order operators, 
Pacific J. Math. {\bf 136} (1989), 355--378.

\bibitem{VG}
A. M. Vershik and V. Ya. Gershkovich,
Nonholonomic dynamical systems, geometry of distributions, and variational 
problems, in : Dynamical Systems VII, Encyclopaedia of Mathematical Sciences
{\sl 16}, Springer-Verlag: Berlin-Heidelberg-New York, 1994, 1--81.

\bibitem{VV}
A. P. Veselov and L. E. Veselova,
Integrable nonholonomic systems on Lie groups,
Math. Notes {\bf 44} (1988), 810--819.

\bibitem{VoG}
S. K. Vodopjanov and A. V. Greshnov,
On extension of functions with bounded mean oscillation onto a space of 
homogeneous type with intrinsic metric,
Siberian Math. Journal {\bf 37} (1996).

\end{thebibliography}
\end{document}